\documentclass[10pt,reqno]{amsart}
\usepackage{geometry} 
\usepackage{graphicx}
\usepackage{amsmath}
\usepackage{color}

\newtheorem{thm}[equation]{Theorem}
\newtheorem{pro}[equation]{Proposition}

\newtheorem{lem}[equation]{Lemma}

\theoremstyle{definition}

\newtheorem{DEF}[equation]{Definition}

\usepackage{amsfonts,amssymb,amscd,amsmath,enumerate,verbatim,calc}

\def\a{\alpha}

\def\e{\epsilon}

\def\LL{\mathcal{L}}

\def\lam{\lambda}

\def\qed{\hfill$\Box$\vspace{3mm}}

\def\ad{\hbox{ad}}

\def\hh{\mathcal H}

\def\bbbz{{\mathbb Z}}

\def\bbbk{{\mathbb K}}

\def\ta{\tilde{\a}}
\def\lam{\Lambda_I^{\bar{i}}}

\def\e{\epsilon}

\def\a{\alpha}
\def\lam{\lambda}

\def\L{L}
\def\dd{\mathcal D}

\begin{document}

\markboth{Double derivations of $n$-Hom-Lie color algebras} {V. Khalili}

\date{}

\centerline{\bf Double derivations of $n$-Hom-Lie color algebras }

\vspace{.5cm}\centerline{Valiollah Khalili\footnote[1]{Department
of mathematics, Faculty of sciences, Arak University, Arak 385156-8-8349, Po.Box: 879, Iran. 
 V-Khalili@araku.ac.ir\\
\hphantom{ppp}2020 Mathematics Subject Classification(s): 17B40, 17B61, 17B75.
\\
\hphantom{ppp}Keywords: $n$-Hom-Lie color algebras, derivations, inner derivations, double derivations }\;\;}

\vspace{1cm} \noindent ABSTRACT:  We study the double derivation algebra $\dd(\LL)$ of $n$-Hom Lie color algebra $\LL$ and  describe the relation between $\dd(\LL)$ and the usual derivation Hom-Lie color algebra $Der(\LL).$  We prove that
the inner derivation algebra $Inn(\LL)$ is an ideal of the double derivation
algebra $\dd(\LL).$ We also show that if $\LL$ is a perfect $n$-Hom Lie color algebra with certain constraints on the base field, then the centralizer of $Inn(\LL)$ in $\dd(\LL)$
is trivial.  In addition, we obtain that for every centerless
perfect $n$-Hom Lie color algebra $\LL,$  the triple derivations of the
derivation algebra $Der(\LL)$ are exactly the derivations of $Der(\LL).$

\vspace{1cm} \setcounter{section}{0}
\section{Introduction}\label{introduction}

The first motivation to study $n$-ary algebras appeared in Physics when
Nambu suggested in 1973 a generalization of Hamiltonian Mechanics with more than one
hamiltonian \cite{N}. The mathematical algebraic foundations of Nambu mechanics have been developed by Takhtajan  in \cite{T1, T2}. The abstract definition of $n$-Lie algebra is due to Filippov in 1985 \cite{F}, then Kasymov \cite{K} deeper investigated their properties. This approach uses the interpretation of Jacobi identity expressing the fact that the adjoint map is a derivation. Nambu mechanics, $n$-Lie algebras revealed to have many applications in physics.

The theory of Hom-algebras originated from Hom-Lie algebras introduced by J. T. Hartwig, D. Larsson and
S. D. Silvestrov in \cite{HLS} in the study of quasi-deformations of Lie algebras of vector fields, including q-deformations
of Witt algebras and Virasoro algebras.
Hom-type generalization of $n$-ary algebras, such as $n$-Hom-Lie algebras and other $n$-ary Hom algebras of Lie type and associative type, were introduced in \cite{AMS1},  by twisting the
identities defining them using an algebra morphism. 
Further properties, construction methods, examples, cohomology and central extensions of $n$-ary Hom-algebras
have been considered in \cite{AMS2, AMS3,  AMM, Y}.

The concept of derivations appear in different mathematical fields with many
different forms. In algebra systems, derivations are linear maps that satisfy the Leibniz relation. There are several kinds of derivations in the theory of Lie algebras, such as Leibniz derivations \cite{KP1},  Jordan derivations
\cite{H} and $\delta-$derivations \cite{F1, F2, F3, K1}. The notion of a generalized derivation  of a Lie algebra and their subalgebras is a generalization of $\delta-$derivation was due to Leger and Luks \cite{LL}. Their results were generalized by many authors. For example, in the case of  color Lie algebras and  Hom-Lie color algebras \cite{CMN, ZF1}, Hom and BiHom–Poisson color algebras \cite{Kh1, Kh2} Lie triple systems and Hom-Lie
triple systems \cite{ZCM1, ZCM2},  color $n$-ary $\Omega$-algebras  and multiplicative $n$-ary Hom-$\Omega-$algebras \cite{KP2, BKP} and many
other works. Another generalization of derivations of Lie algebras are triple derivations and generalized  triple derivations. It was first
introduced independently  by Muller where it was called prederivation \cite{Mu}. It can be easily checked that, for any Lie algebra, every derivation is a  triple derivation, but the converse does not always hold \cite{Z}. The author proved that every triple derivation of a perfect Lie algebra with zero
center is a derivation and every derivation of the derivation algebra is an
inner derivation. In recent years,
much progress has been obtained on the Lie triple derivations and $n$-Lie drivations of Lie algebras (see \cite{AW, LC, LW, M, WX, XW, Z, ZCM, ZWC}). 

The double derivations provide a useful tool to study $n$-Lie algebras by making use of linear mappings. Double derivations are similar to the triple derivations of Lie algebras to some extent \cite{Z}. Recently, in \cite{BGZL} and  \cite{SCZ}  the authors studied double derivations of $n$-Lie algebras and $n$-Lie superalgebras. The main purpose of our work is to consider the double derivations of multiplicative $n$-Hom- Lie color algebras.
 
The paper is organized as follows. Section 2 contains some necessary  basic definitions and notions. In this section,  we show that the double derivation algebra of a multiplicative  $n$-Hom-Lie color algebra $\LL$ is a Hom-Lie coloralgebra and
for a perfect $n$-Hom-Lie color algebra $\LL$ its inner derivation algebra is an ideal of the double
derivation algebra. In Section 3, we construct a new double derivation of a perfect
centerless multiplicative $n$-Hom-Lie color algebra $\LL$  and obtain a Hom-Lie color algebra homomorphism.
In section 4, we study triple derivations of the derivation algebra $Der(\LL)$ and the inner
derivation algebra $Inn(\LL).$

\section{Double Derivations of  multiplicative $n$-Hom-Lie color algebras} \setcounter{equation}{0}\

Let us begin with some necessary important basic definitions and notations on graded spaces, 
graded algebras and  $n$-Hom-Lie color algebras used in other sections.  For a detailed discussion of this subject, we refer the reader to the literature \cite{BS}. We also   prove that
the inner derivation algebra $Inn(\LL)$ of $\LL$ is an ideal of the double derivation
algebra of $\LL.$

Let $G$ be any additive abelian group,  a vector space $V$ is
said to be {\em $G$-graded}, if there is a family
$\{V_g\}_{g\in G}$ vector subspaces such that
$V=\bigoplus_{g\in G}V_g.$ An element $v\in V$ is said to
be {\em homogeneous of the degree $g$} if $v\in
V_g,~g\in G,$ and in this case, $g$ is called the {\em
	color} of $v.$ We denote
$\hh(V)$ the set of all homogeneous elements in $V.$ 

Let $V=\bigoplus_{g\in G}V_g$ and
$W=\bigoplus_{h\in G}W_h$ be two $G-$graded vector
spaces. A linear mapping $f : V\longrightarrow W$ is said to be
{\em homogeneous of degree} $\theta\in G$ if
$$
f(V_g)\subset W_{g+\theta}, ~~~\forall g\in G.
$$
If in addition, $f$ is homogeneous of degree zero, namely,
$f(V_g)\subset W_g$ holds for any $g\in G,$ then we call $f$
is  {\em even}.

An algebra $A$ (with the juxtaposition product) is said to be $G$-graded if its underlying
vector space is $G$-graded, i.e.,
$A=\bigoplus_{g\in G}A_g,$ and if $A_g A_h\subset
A_{g+h},$ for $g, h\in G.$ A subalgebra of $A$ is said
to be graded if it is a graded as a subspace of $A.$

Let $B$ be another $G$-graded algebra. A morphism $\varphi
: A\longrightarrow B$ of $G$-graded algebras is a homomorphism
of the algebra $A$ into the algebra $B,$ which is an even mapping.

\begin{DEF}\label{bi} Let $G$ be an addtive abelian group. A map $\e :
	G\times G\longrightarrow\bbbk\setminus\{0\}$ is called a
	skew-symmetric {\em bi-character} on $G$ if for all $g,h,
	k\in G,$
	\begin{itemize}
		\item[(i)] $\e(g, h)\e(h, g)=1,$
		
		\item[(ii)] $\e(g+h, k)=\e(g, k) \e(h, k),$
		
		\item[(iii)]  $\e(g, h+k)=\e(g, h)\e(g, k).$
	\end{itemize}
\end{DEF}
The definition above implies that in particular, the following
relations hold
$$
\e(g, 0)=1=\e(0, g),~~~~\e(g,
g)=1(\hbox{or}~-1),~~~\forall g\in G.
$$

Throughout this paper, if $x$ and $ y$  are two homogeneous elements
of a $G-$graded vector space and $|x|, |y|$  
are  their degrees respectively, then for
convenience, we write $\e(x, y)$ instead of $\e(|x|, |y|).$ It is worth mentioning that,  We
unless otherwise stated, in the sequel all the graded spaces are over
the same abelian group $G$ and the bi-character will be the
same for all structures.

\begin{DEF}\cite{F} An {\em $n$-Lie algebra} is a linear spaces $\L$ equipped with $n$-ary operation  satisfies the following identity:
\begin{itemize}
\item[(i)]  $[x_1, ... x_i, x_{i+1}, ..., x_n]=-[x_1, ... x_{i+1}, x_{i}, ..., x_n], ~~i=1, ..., n,$
\item[(ii)]	$
	\big[x_1,...,x_{n-1},[y_1,...,y_n]\big] =\sum_{i=1}^n\big[y_1,...,y_{i-1},[x_1,...,x_{n-1}, y_i
	], y_{i+1},... y_n\big],
	$
\end{itemize}
for all  $x_i, y_j\in\hh(\LL).$
\end{DEF}

\begin{DEF}\cite{BS} An {\em $n$-Hom-Lie color algebra} is a graded linear space $\LL=\bigoplus_{g\in G}\LL_g$ with
	an $n$-linear map $[., ..., .]~:~  \LL\times...\times\LL\longrightarrow\LL,$ a bicharacter $\e :
	G\times G\longrightarrow\bbbk\setminus\{0\}$ and an even linear
	map 	$\a : \LL\longrightarrow\LL$ such that  
	\begin{itemize}
		\item[(i)] $[x_1, ... x_i, x_{i+1}, ..., x_n]=-\e(x_i, x_{i+1})[x_1, ... x_{i+1}, x_{i}, ..., x_n], ~~i=1, ..., n,$
		
		\item[(ii)]
		\begin{eqnarray}
			&\big[\a(x_1),...,\a(x_{n-1}),[y_1, y_2, ..., y_n]\big] =\\
			\nonumber&\sum_{i=1}^n\e(X, Y_i)\big[\a(y_1), ..., \a(y_{i-1}), [x_1, ..., x_{n-1}, y_i
			], \a(y_{i+1}), ..., \a(y_n)\big],
		\end{eqnarray}		
	\end{itemize}
	where $x_i, y_j\in\hh(\LL), X =\sum_{i=1}^{n-1} x_i, Y_i =\sum_{j=1}^{i}y_{j-1}$ and $y_0 =e.$
\end{DEF}

An $n$-Hom-Lie color algebra $(\LL, [., ..., .], \a, \e)$ is said to be multiplicative if $\a$ is an endomorphism, i.e. a linear map on $\LL$ which is also a homomorphism with
respect to multiplication $[., ..., .].$ We also call it a regular Hom-Lie color algebra  if $\a$ is an automorphism. We recover $n$-Lie 
algebra when we have $\a=id$ and $G=\{e\}.$ If $\a = id_{\LL}$ and  $G=\bbbz_2$ with $\e(x, y)=(-1)^{x y}$ for any $ x, y \in\hh(\LL),$ we get $n$-Lie super algebra.  For some standard exampls, we refer the reader can be found in \cite{BS}.

\begin{DEF} 
\begin{itemize}
\item[(1)] A graded subspace $S$ of an $n$-Hom-Lie color algebra $\LL$ is  a   $n$-Hom-Lie color subalgebra of $\LL$ if $\a(S)\subseteq S$ and $ [S, S, ..., S]\subseteq S.$
\item[(2)]	An $n$-Hom-Lie  color ideal $I$ of an $n$-Hom-Lie color algebra $\LL$ is a graded subspace of $\LL$ such that $\a(I)\subseteq I$ and $ [I, \LL, ..., \LL]\subseteq I.$
\end{itemize}
\end{DEF}

Let $S_1, S_2, ..., S_n$ be $n$-Hom-Lie subalgebras of an $n$-Hom-Lie color algebra $\LL.$
Denote by $[S_1, S_2, ..., S_n]$ the  $n$-Hom-Lie subalgebra of $\LL$ generated by all elements $[x_1, x_2, ..., x_n],$ where $x_i \in S_i, i = 1, 2, ..., n.$  The algebra $[\LL,  . . . , \LL ]$ is called the derived algebra of $\LL,$ and is denoted by $\LL^1.$ If $\LL^1=\LL,$ then $\LL$ is called a perfect $n$-Hom-Lie color algebra. The center of an $n$-Hom-Lie color algebra $\LL$ is denoted by $Z(\LL) = \{x \in\LL~:~  [x, \LL,  . . . , \LL ] = 0\}$ which is a $n$-Hom-Lie color ideal of $\LL$ (see Theorem 2.16 in \cite{BS}). For a subset $S$ of $\LL,$
the centralizer of $S$ in $\LL$ is defined by $C_{\LL}(S) = \{x \in\LL~:~[x, S, \LL, . . . , \LL] = 0\}.$ Note that $C_{\LL}(\LL)=Z(\LL).$ 

From now on, we let   $(\LL, [., ..., .], \a, \e)$  be a multiplicative $n$-Hom-Lie color algebra, which we will denote by $\LL,$ for short. For
any non-negative integer $k,$ denote by $\a^k$  the $k-$times
composition of $\a.$  In particular $\a^0=id, ~\a^1=\a.$

In the  following we give a more general definition of derivations  and related objects.

\begin{DEF}\cite{BS} For any $k\geq 1,$ we call $D\in End(\LL)$ an {\em $\a^k$-derivation} of degree $d$ of the multiplicative $n$-Hom-Lie color algebra $\LL$ if
\begin{itemize}
\item[(1)]  $D\circ\a=\a\circ D.$ 
\item[(2)]  For all $x_1, x_2, ..., x_n\in\hh(\LL),$
\begin{equation*}	
D\big([x_1, ..., x_n]\big)=\sum_{i=1}^{n}\e(d, X_i)\big[\a^k(x_1), ..., \a^k(x_{i-1}),   D(x_i), \a^k(x_{i+1}), ..., \a^k(x_n)\big],
\end{equation*}
where $ X_i =\sum_{j=1}^{i} x_j.$
\end{itemize}	
\end{DEF}

We denote the set of all $\a^k$-derivations of the multiplicative $n$-Hom-Lie color algebra $\LL$ by $Der_{\a^k}(\LL).$  For any $D\in Der_{\a^k}(\LL)$  and  $D'\in Der_{\a^s}(\LL),$ let us define their color commutator $[D, D']=D\circ D'-\e(d, d')D'\circ D.$ Then $[D, D']\in Der_{\a^{k+s}}$ (see Lemma 5.6 in \cite{BS}). Denote by $Der(\LL) = \bigoplus_{k\geq -1}Der_{\a^k} (\LL).$ We also have that $(Der(\LL), [., .], \ta=D\circ\a, \e)$ is a Hom-Lie color algebra (see Proposition 5.7 in \cite{BS}).

For $x_1, ..., x_{n-1}\in\hh(\LL)$ satisfying $\a(x_i)=x_i, i=1, 2, ..., n-1,$ we define the map $\ad_k(x_1, ..., x_{n-1})~:~\LL\longrightarrow\LL$  by
$$
\ad_k(x_1, ..., x_{n-1})(y) :=[x_1, ..., x_{n-1}, \a^k(y)], ~\hbox{~and~}~k\ge 1,
$$
for all $y\in\hh(\LL).$ Then $\ad_k(x_1, ..., x_{n-1}),$ is an $\a^{k+1}$-derivation of $\LL$ (see Lemma 2.2 in \cite{AMM}).   We call $\ad_{k}(x_1, ..., x_{n-1})$ an {\em inner $\a^{k+1}$-derivation}. Denote by  $Inn_{\a^k}(\LL)$ the space generate by all the inner $\a^{k+1}$-derivations.  Set  $Inn(\LL)=\bigoplus_{k\geq 0}Inn_{\a^k}(\LL).$  It is not difficult to show that $(Inn(\LL), [., .], \ta=D\circ\a, \e)$  is a Hom-Lie color algebra. In addition, $Inn(\LL)$ is a Hom-Lie color ideal of $Der(\LL)$ (see Proposition 2.4  in \cite{AMM}).
\begin{DEF}	Let   $(\LL, [., ..., .], \a, \e)$  be a multiplicative $n$-Hom-Lie color algebra with $n\geq 3.$ A linear
map $D~:~ \LL\longrightarrow\LL$ is called a {\em double $\a^k$-derivation} of degree $d$ of $\LL$ if 
\begin{itemize}
\item[(1)]  $D\circ\a=\a\circ D,$ 
\item[(2)]  for all $x_i, y_j\in\hh(\LL),~i=1, 2, ..., n-1,~j=1, 2, ..., n,$
\begin{eqnarray*}
&~&D\big(\big[x_1, ..., x_{n-1}, [y_1, ..., y_{n}]\big]\big)\\
&=&\sum_{i=1}^{n-1}\e(d, X_i)\big[\a^k(x_1), ..., D(x_i), ..., \a^k(x_{n-1}), [\a^k(y_1), ..., \a^k(y_n)]\big]\\
&+&\sum_{j=1}^{n}\e(d, X+Y_{j})\big[\a^k(x_1), ..., \a^k(x_{n-1}),  [\a^k(y_1), ..., D(y_j), ..., \a^k(y_n)]\big].
\end{eqnarray*}
\end{itemize}	
\end{DEF}
Let us denote the set of all double $\a^k$-derivations of the multiplicative $n$-Hom-Lie color algebra $\LL$ by $\dd_{\a^k}(\LL)$ and denote by $\dd(\LL) = \bigoplus_{k\geq -1}\dd_{\a^k} (\LL),$  the vector space spanned by the double derivations of $\LL.$ We call it {\em  double derivation algebra} of $\LL$ (see Theorem \ref{Doublealg} below).	

It is obvious that the derivations of an $n$-Hom-Lie color algebra $\LL$ are double derivations, and hence we have $Der(\LL)\subset\dd(\LL).$ Note that there exists an $n$-Lie superalgebra $\LL$ such that $Der(\LL)\neq\dd(\LL)$ (see Remark 2.5 in \cite{SCZ}).

\begin{thm}\label{Doublealg}	Let   $(\LL, [., ..., .], \a, \e)$  be a multiplicative $n$-Hom-Lie color algebra. Then $\dd(\LL)$ is a  Hom-Lie color subalgebra of  the general linear Hom-Lie color algebra $(gl(\LL), [., .], \ta, \e)$ with $\ta(D)=D\circ\a$ for all $D\in\dd(\LL).$	
\end{thm}
\noindent {\bf Proof.} First, let us prove that $\ta(\dd)(\LL)\subseteq\dd(\LL).$ Suppose that $D\in\dd_{\a^k}(\LL),$ for any $x_1, x_2, ..., x_{n-1}, y_1, ..., y_n\in\hh(\LL),$ we have
\begin{eqnarray*}
&\ta(D)\bigg(\bigg[x_1, ..., x_{n-1}, [y_1, ..., y_{n}]\bigg]\bigg)\\
=&(D\circ\a)\bigg(\bigg[x_1, ..., x_{n-1}, [y_1, ..., y_{n}]\bigg]\bigg)\\
=&D\bigg(\big[\a(x_1), ..., \a(x_{n-1}), [\a(y_1), ..., \a(y_{n})]\big]\bigg)\\
=&\sum_{i=1}^{n-1}\e(d, X_i)\\
&.\big[\a^{k+1}(x_1), ..., D(\a(x_i)), ..., \a^{k+1}(x_{n-1}), [\a^{k+1}(y_1), ..., \a^{k+1}(y_n)]\big]\\
&+\sum_{j=1}^{n}\e(d, X+Y_j)\\
&.\big[\a^{k+1}(x_1), ..., \a^{k+1}(x_{n-1}), [\a^{k+1}(y_1), ..., D(\a(y_j)), ..., \a^{k+1}(y_n)]\big]\\
=&\sum_{i=1}^{n-1}\e(d, X_i)\\
&.\big[\a^{k+1}(x_1), ..., (D\circ\a)(x_i), ..., \a^{k+1}(x_{n-1}), [\a^{k+1}(y_1), ..., \a^{k+1}(y_n)]\big]\\
&+\sum_{j=1}^{n}\e(d, X+Y_j)\\
&.\big[\a^{k+1}(x_1), ..., \a^{k+1}(x_{n-1}), [\a^{k+1}(y_1), ..., (D\circ\a)(y_j), ..., \a^{k+1}(y_n)]\big]\\
=&\sum_{i=1}^{n-1}\e(d, X_i)\\
&.\big[\a^{k+1}(x_1), ..., (\ta(D))(x_i), ..., \a^{k+1}(x_{n-1}), [\a^{k+1}(y_1), ..., \a^{k+1}(y_n)]\big]
\end{eqnarray*}
\begin{eqnarray*}
&+\sum_{j=1}^{n}\e(d, X+Y_j)\\
&.\big[\a^{k+1}(x_1), ..., \a^{k+1}(x_{n-1}), [\a^{k+1}(y_1), ..., (\ta(D))(y_j), ..., \a^{k+1}(y_n)]\big].	
\end{eqnarray*}
This means that $\ta(D)$ is a duble $\a^{k+1}$-derivation i.e. $\ta(D)\in\dd_{\a^{k+1}}(\LL).$

Next, assum that $D_{1}\in\dd_{\a^k}(\LL)$ of degree $d_1$ and $D_{2}\in\dd_{\a^s}(\LL)$ of degree $d_2.$ It is clear that $[D_1, D_2]\circ\a=\a\circ[D_1, D_2].$ For any $x_1, x_2, ..., x_{n-1}, y_1, ..., y_n\in\hh(\LL),$ we  also have 
\begin{eqnarray*}
& D_{1}D_{2}\bigg(\big[x_1, ..., x_{n-1}, [y_1, ..., y_{n}]\big]\bigg)\\
=& D_1\bigg(\sum_{i=1}^{n-1}\e(d_2, X_i)\big[\a^s(x_1), ..., D_2(x_i), ..., \a^s(x_{n-1}), [\a^s(y_1), ..., \a^s(y_n)]\big]\\
&+\sum_{j=1}^{n}\e(d_2, X+Y_j)\\
&.\big[\a^{s}(x_1), ..., \a^s(x_{n-1}), [\a^s(y_1), ..., D_2(y_j), ..., \a^s(y_n)]\big]\bigg)\\
=&\sum_{i=1}^{n-1}\e(d_1+d_2, X_i)\\
&.\big[\a^{k+s}(x_1), ..., D_1D_2(x_i), ..., \a^{k+s}(x_{n-1}), [\a^{k+s}(y_1), ...,  \a^{k+s}(y_n)]\big]\\
&+\sum_{i=1}^{n-1}\sum_{r<i}\e(d_1, X_r)\e(d_2, X_i)\\
&.\bigg[\a^{k+s}(x_1), ..., D_1(\a^{s}(x_r)), ..., D_2(\a^k(x_i)), ..., \a^{k+s}(x_{n-1}),\\
&~\big[\a^{k+s}(y_1), ..., \a^{k+s}(y_n)\big]\bigg]\\
&+\sum_{i=1}^{n-1}\sum_{i<r}\e(d_1, X_r)\e(d_2, X_i)\e(d_1, d_2)\\
&.\bigg[\a^{k+s}(x_1), ..., D_2(\a^{k}(x_i)), ...,  D_1(\a^s(x_r)), ..., \a^{k+s}(x_{n-1}),\\
&~\big [\a^{k+s}(y_1), ..., \a^{k+s}(y_n)\big]\bigg]\\
&+\sum_{i=1}^{n-1}\sum_{t=1}^{n-1}\e(d_2, X_i)\e(d_1, X)\e(d_1, d_2)\\
&.\bigg[\a^{k+s}(x_1), ..., D_2(\a^k(x_i)),  ..., \a^{k+s}(x_{n-1}),\\
&~\big[\a^{k+s}(y_1), ..., D_1(\a^s(y_t)), ..., \a^{k+s}(y_n)\big]\bigg]\\
&+\sum_{j=1}^{n}\sum_{m=1}^{n-1}\e(d_1, X_m)\e(d_2, X+Y_j)\\
&.\bigg[\a^{k+s}(x_1), ..., D_1(\a^s(x_j)), ...,\a^{k+s}(x_{n-1}),\\
&~\big[\a^{k+s}(y_1), ..., D_2(\a^k(y_m)), ..., \a^{k+s}(y_n)\big]\bigg]
\end{eqnarray*}
\begin{eqnarray*}
&+\sum_{j=1}^{n}\sum_{l< j}\e(d_1, X+Y_l)\e(d_2, X+Y_j)\\
&.\bigg[\a^{k+s}(x_1), ..., \a^{k+s}(x_{n-1}),\\
&~\big[\a^{k+s}(y_1), ..., D_1(\a^s(y_l)), ..., D_2(\a^k(y_j)), ..., \a^{k+s}(y_n)\big]\bigg]\\
&+\sum_{j=1}^{n}\sum_{j<l}\e(d_1, X+Y_l)\e(d_2, X+Y_j)\e(d_1, d_2)\\
&.\bigg[\a^{k+s}(x_1), ..., \a^{k+s}(x_{n-1}),\\
&~\big[\a^{k+s}(y_1), ..., D_2(\a^k(y_j)), ..., D_1(\a^s(y_l)), ..., \a^{k+s}(y_n)\big]\bigg]\\
&+\sum_{j=1}^{n}\e(d_1+d_2, X+Y_j)\\
&.\bigg[\a^{k+s}(x_1), ..., \a^{k+s}(x_{n-1}), \big[\a^{k+s}(y_1),  ..., D_1D_2(x_j), ..., \a^{k+s}(y_n)\big]\bigg].
\end{eqnarray*}
Similarly,
\begin{eqnarray*}
& \e(d_1, d_2)D_{2}D_{1}\bigg(\big[x_1, ..., x_{n-1}, [y_1, ..., y_{n}]\big]\bigg)\\
=&\sum_{i=1}^{n-1}\e(d_1+d_2, X_i)\\
&.\bigg[\a^{k+s}(x_1), ..., D_2D_1(x_i), ..., \a^{k+s}(x_{n-1}), \big[\a^{k+s}(y_1),  ..., \a^{k+s}(y_n)\big]\bigg]\\
&+\sum_{i=1}^{n-1}\sum_{r<i}\e(d_2, X_r)\e(d_1, X_i)\\
&.\bigg[\a^{k+s}(x_1), ..., D_2(\a^{k}(x_r)), ..., D_1(\a^s(x_i)),  ..., \a^{k+s}(x_{n-1}),\\
~&\big[\a^{k+s}(y_1), ..., \a^{k+s}(y_n)\big]\bigg]\\
&+\sum_{i=1}^{n-1}\sum_{i<r}\e(d_2, X_r)\e(d_1, X_i)\e(d_2, d_1)\\
&.\bigg[\a^{k+s}(x_1), ..., D_1(\a^{s}(x_i)),  ..., D_2(\a^k(x_r)), ..., \a^{k+s}(x_{n-1}),\\
~&\big[\a^{k+s}(y_1), ..., \a^{k+s}(y_n)\big]\bigg]\\
&+\sum_{i=1}^{n-1}\sum_{t=1}^{n-1}\e(d_1, X_i)\e(d_2, X)\e(d_2, d_1)\\
&.\bigg[\a^{k+s}(x_1), ..., D_1(\a^s(x_i)),  ..., \a^{k+s}(x_{n-1}),\\
~&\big [\a^{k+s}(y_1), ..., D_2(\a^k(y_t)), ..., \a^{k+s}(y_n)\big]\bigg]
\end{eqnarray*}
\begin{eqnarray*}
&+\sum_{j=1}^{n}\sum_{m=1}^{n-1}\e(d_2, X_m)\e(d_1, X+Y_j)\\
&.\bigg[\a^{k+s}(x_1), ..., D_2(\a^k(x_j)),  ..., \a^{k+s}(x_{n-1}),\\
~&\big[\a^{k+s}(y_1), ..., D_1(\a^k(y_m)), ..., \a^{k+s}(y_n)\big]\bigg]\\
&+\sum_{j=1}^{n}\sum_{l< j}\e(d_2, X+Y_l)\e(d_1, X+Y_j)\\
&.\bigg[\a^{k+s}(x_1), ..., \a^{k+s}(x_{n-1}),\\
~&\big[\a^{k+s}(y_1), ..., D_2(\a^k(y_l)), ..., D_1(\a^s(y_j)), ..., \a^{k+s}(y_n)\big]\bigg]\\
&+\sum_{j=1}^{n}\sum_{j<l}\e(d_2, X+Y_l)\e(d_1, X+Y_j)\e(d_2, d_1)\\
&.\bigg[\a^{k+s}(x_1), ..., \a^{k+s}(x_{n-1}),\\ 
~&\big[\a^{k+s}(y_1), ..., D_1(\a^k(y_j)), ..., D_2(\a^k(y_l)), ..., \a^{k+s}(y_n)\big]\bigg]\\
&+\sum_{j=1}^{n}\e(d_1+d_2, X+Y_j)\\
&.\bigg[\a^{k+s}(x_1), ..., \a^{k+s}(x_{n-1}), \big[\a^{k+s}(y_1), ...,  D_2D_1(x_j), ..., \a^{k+s}(y_n)\big]\bigg].
\end{eqnarray*}
Then from easy computation we have
\begin{eqnarray*}
&[D_{1}, D_{2}]\bigg(\big[x_1, ..., x_{n-1}, [y_1, ..., y_{n}]\big]\bigg)\\
=&\big( D_1  D_2-\e(d_1, d_2)D_{2}D_{1}\big)\bigg(\big[x_1, ..., x_{n-1}, [y_1, ..., y_{n}]\big]\bigg)\\
=& \sum_{i=1}^{n-1}\bigg[\a^s(x_1), ..., \big(D_1  D_2-\e(d_1, d_2)D_{2}D_{1}\big)(x_i), ..., \a^s(x_{n-1}),\\
~&\big[\a^s(y_1), ..., \a^s(y_n)\big]\bigg]\\
&+\sum_{j=1}^{n}\e(d_1+d_2, X+Y_j)\\
&.\bigg[\a^{k+s}(x_1), ..., \a^{k+s}(x_{n-1}),\\
~&\big[\a^{k+s}(y_1), ...,  \big(D_1  D_2-\e(d_1, d_2)D_{2}D_{1}\big)(y_j), ..., \a^{k+s}(y_n)\big]\bigg]\\
=&\sum_{i=1}^{n-1}\bigg[\a^s(x_1), ..., \big[D_1,  D_2\big](x_i), ..., \a^s(x_{n-1}), \big[\a^s(y_1), ..., \a^s(y_n)\big]\bigg]\\ 
&+\sum_{j=1}^{n}\e(d_1+d_2, X+Y_j)\\
&.\bigg[\a^{k+s}(x_1), ..., \a^{k+s}(x_{n-1}), \big[\a^{k+s}(y_1), ..., \big[D_1,  D_2\big](x_j), ..., \a^{k+s}(y_n)\big]\bigg].
\end{eqnarray*}
This implies that, $[D_1, D_2]$ is a  double $\a^{k+s}$-derivation of degree $d_1+d_2$ and so $[D_1, D_2]\in \dd(\LL)$ as required.\qed

\begin{thm}\label{ideal} If $(\LL, [., ..., .], \a, \e)$ is a perfect  multiplicative $n$-Hom-Lie color algebra, then the inner derivation algebra of $\LL$ is a Hom-Lie color ideal of the Hom-Lie color algebra $(\dd(\LL), [., .], \ta, \e).$ 
\end{thm}	
\noindent {\bf Proof.}  First, we are going to check thatt  $\ta(Inn(\LL))\subseteq Inn(\LL).$ For $x_1, x_2, ..., x_{n-1}\in\hh(\LL)$ and any $y\in\hh(\LL),$  we have  
\begin{eqnarray}
\nonumber\ta\bigg(\ad_{k}\big(x_1, ..., x_{n-1}\big)\bigg)(y)&=&\ad_{k}\big(x_1, ..., x_{n-1}\big)\big(\a(y)\big)\\
\nonumber&=&\big[x_1, ..., x_{n-1}, \a^{k+1}(y)\big]\\
\nonumber&=&\ad_{k+1}\big(x_1, ..., x_{n-1}\big)(y).	
\end{eqnarray}
This means that $\ta\big(\ad_{k}(x_1, ..., x_{n-1})\big)$ is an inner $\a^{k+2}$-derivation i.e. $\ta\big(Inn(\LL)\big)\subseteq Inn(\LL).$

Next, let $D\in\dd_{\a^s}(\LL)$ and $ad_k(x_1, x_2, ..., x_{n-1})$ be an inner $\a^{k+1}$-derivation. Since $\LL$ is perfect, there exists a finite index set $I$ and $x_{i_j}\in\LL, i \in I, j=1, ..., n,$ such that
$$
x_1=\sum_{i\in I}[x_{i_1}, x_{i_2}, ..., x_{i_n}].
$$
For any $y\in\LL$ we have
\begin{eqnarray*}
&\big[D, ad_k(x_1, ..., x_{n-1})\big](y)\\
=&D\circ ad_{k}(x_1, ..., x_{n-1})(y)-\e(d, X)ad_{k}(x_1, ..., x_{n-1})\circ D(y)\\
=&D\big([\a^k(x_1), ..., \a^k(x_{n-1}), \a^k(y)]\big)-\e(d, X)\big[\a^k(x_1), ..., \a^k(x_{n-1}), \a^k(D(y))\big]\\
=&\sum_{i\in I}D\bigg(\big[[\a^k(x_{i_1}), \a^k(x_{i_2}), ..., \a^k(x_{i_n})], \a^k(x_2), ..., \a^k(x_{n-1}), \a^k(y)\big]\bigg)\\ 
&-\e(d, X)\big[\a^{k+s}(x_1), ..., \a^{k+s}(x_{n-1}), D(\a^k(y))\big]\\
=&\sum_{i\in I}\bigg(\sum_{j=2}^{n-2}\e(d, x_2+...+x_{j-1})\\
&.\big[\a^{k+s}(x_2), ..., D(x_j),..., \a^{k+s}(x_{n-1}), \a^{k+s}(y), [\a^{k+s}(x_{i_1}), ..., \a^{k+s}(x_{i_n})]\big]\\
&+\e(d, x_2+...x_{n-1})\\
&.\big[\a^{k+s}(x_2), ..., \a^{k+s}(x_{n-1}),  D(\a^k(y)), [\a^{k+s}(x_{i_1}), ..., \a^{k+s}(x_{i_n})]\big]\\
&+\sum_{t=1}^{n}\e(d, x_2+...+x_{n-1}+y+X_{i_{t}})\\
&.\big[\a^{k+s}(x_2),  ..., \a^{k+s}(x_{n-1}), \a^{k+s}(y), [\a^{k+s}(x_{i_1}), ..., D(x_{i_t}), ..., \a^{k+s}(x_{i_n})]\big]\bigg)\\
&-\e(d, X)\big[\a^{k+s}(x_1), ..., \a^{k+s}(x_{n-1}), D(\a^k(y))\big]
\end{eqnarray*}
\begin{eqnarray*}
=&\sum_{i\in I}\bigg(\sum_{j=2}^{n-1}\e(d, X_j)\\
&.\big[[\a^{k+s}(x_{i_1}), ..., \a^{k+s}(x_{i_n})], \a^{k+s}(x_2), ..., D(x_j), ..., \a^{k+s}(x_{n-1}), \a^{k+s}(y)\big]\\
&+\e(d, X)\big[[\a^{k+s}(x_{i_1}), ..., \a^{k+s}(x_{i_n})], \a^{k+s}(x_2), ..., \a^{k+s}(x_{n-1}),  D(\a^k(y))\big]\bigg)\\
&+\sum_{i\in I}\bigg(\sum_{t=1}^{n}\e(d, X_{i_t})\\
&.\big[[\a^{k+s}(x_{i_1}), ..., \a^{k+s}(x_{i_n})], \a^{k+s}(x_2), ..., \a^{k+s}(x_{n-1}),\a^{k+s}(y)\big]\bigg)\\
&-\e(d, X)\big[\a^{k+s}(x_1), ..., \a^{k+s}(x_{n-1}), D(\a^k(y))\big]\\
=&\sum_{j=2}^{n-1}\e(d, X_j)\big[\a^{k+s}(x_1), ..., D(x_j), ..., \a^{k+s}(x_{n-1}), \a^{k+s}(y)\big]\\
&+\sum_{i\in I}\bigg(\sum_{t=1}^{n}\e(d, X_{i_t})\\
&.\big[[\a^{k+s}(x_{i_1}), ..., D(x_{i_t}), ..., \a^{k+s}(x_{i_n})], \a^{k+s}(x_2), ..., \a^{k+s}(x_{n-1}), \a^{k+s}(y)\big]\bigg)\\
=&\sum_{j=2}^{n-1}\e(d, X_j)ad_{k+s}\big(x_1, ..., D(x_j), ..., x_{n-1}\big)(y)\\
&+\sum_{i\in I}\sum_{t=1}^{n}\e(d, X_{i_t})ad_{k+s}\big([x_{i_1}, ..., D(x_{i_t}), ..., x_{i_n}], x_2, ..., x_{n-1}\big)(y).
\end{eqnarray*}

Consequently, $[D, ad_k(x_1, ..., x_{n-1}]\in Inn_{\a^{k+s+1}}(\LL),$ which implies that $Inn(\LL)$ is a Hom-Lie color  ideal of Hom-Lie color algebra $(\dd(\LL), [., .], \ta, \e).$\qed

Note that, if $\LL$ is not a perfect  multiplicative $n$-Hom-Lie color algebra, then $Inn(\LL)$ may not be an ideal of $\dd(\LL)$ (see Rimark 2.8 in \cite{SCZ}).\\

\section{Double Derivations of Perfect multiplicative $n$-Hom-Lie color algebras} \setcounter{equation}{0}\

In this section, we focus on the Perfectnes of   multiplicative $n$-Hom-Lie color algebras by centering our attention in
those of centerless.

Let $(\LL, [., ..., .], \a, \e)$ be an arbitrary centerless multiplicative $n$-Hom-Lie color algebra. For all $D\in\dd_{\a^k}(\LL),$  we define a linear map $\delta_D~:~\LL\longrightarrow\LL,$ with $\delta_D\circ\a=\a\circ\delta_D$ by 
\begin{equation}\label{122}
	\delta_D(x)=\left\{
	\begin{array}{ll}
		D(x)&;~x\in\LL\setminus\LL^1\\
		\sum_{i\in I}\sum_{t=1}^{n}\e(d, X_{i_t})[\a^k(x_{i_1}), ..., D(x_{i_t}), ..., \a^k(x_{i_n})]&;~x=\sum_{i\in I}[x_{i_1},\\& ..., x_{i_n}]\in\LL^1,
	\end{array}
	\right.
\end{equation}
for all $x\in\hh(\LL).$
\begin{lem}
The linear map $\delta_D$  defined as Eq. (\ref{122}) is well defined.
\end{lem}
\noindent {\bf Proof.} Suppose that 
$$
x=\sum_{i\in I}[x_{i_1}, ..., x_{i_n}]=\sum_{j\in J}[x_{j_1}, ..., x_{j_n}]=y,
$$
we set $\delta_D(x)=A$ and $\delta_D(y)=B.$ For $D\in\dd_{\a^k}(\LL)$ of degree $d$ and $z_1, ..., z_{n-1}\in\hh(\LL),$ we have
\begin{eqnarray*}
&~&\big[\a^k(z_1), ..., \a^k(z_{n-1}), A\big]\\
&=&\bigg[\a^k(z_1), ..., \a^k(z_{n-1}), \sum_{i\in I}\sum_{t=1}^{n}\e(d, X_{i_t})[\a^k(x_{i_1}), ..., D(x_{i_t}), ..., \a^k(x_{i_n})]\bigg]\\
&=&\e(d, Z)\bigg(D\big([z_1, ..., z_{n-1}, x]\big)-\\
&~&~~\sum_{s=1}^{n-1}\e(d, Z_s)\big[\a^k(z_1), ..., D(z_s) , ..., \a^k(z_{n-1}), \a^k(x)\big]\bigg)\\
&=&\e(d, Z)\bigg(D\big([z_1, ..., z_{n-1}, y]\big)-\\
&~&~~\sum_{s=1}^{n-1}\e(d, Z_s)\big[\a^k(z_1), ..., D(z_s) , ..., \a^k(z_{n-1}), \a^k(y)\big]\bigg)\\
&=&\e(d, Z)\bigg(D\bigg(\big[z_1, ..., z_{n-1}, \sum_{j\in J}[x_{j_1}, ..., x_{j_n}]\big]\bigg)-\\
& ~&~~\sum_{s=1}^{n-1}\e(d, Z_s)\big[\a^k(z_1), ..., D(z_s) , ..., \a^k(z_{n-1}), \sum_{j\in J}[\a^k(x_{j_1}), ..., \a^k(x_{j_n})]\big]\bigg)\\
&=&\bigg[\a^k(z_1), ..., \a^k(z_{n-1}), \sum_{j\in J}\sum_{s=1}^{n}\e(d, Y_{j_s})\big[\a^k(y_{j_1}), ..., D(y_{j_s}), ..., \a^k(y_{j_n})\big]\bigg]\\
&=&\big[\a^k(z_1), ..., \a^k(z_{n-1}), B\big].
\end{eqnarray*}
Therefore, we get that $[\a^k(z_1), ..., \a^k(z_{n-1}), A-B]=0.$ This implies $A-B\in Z(\LL)=0,$ we obtain $\delta_D(x)=\delta_D(y).$ Thus, $\delta_D$  is well-defined.\qed\\  

Consequently, we obtain a linear map $\delta ~:~\dd(\LL)\longrightarrow End(\LL)$ defined by
$$
\delta(D)=\delta_D, ~~ \forall D\in\dd(\LL).
$$
\begin{pro} Let $(\LL, [., ..., .], \a, \e)$ be a centerless perfect multiplicative $n$-Hom-Lie color algebra. Suppose that $D\in\dd(\LL),$ then the following assertions hold.\\
	\begin{itemize}
		\item[(1)]	 $\delta_D$ is a double derivation of $\LL.$
		
		\item[(2)] For all $x_1, ..., x_n\in\hh(\LL),$ we have 
		\begin{equation*}
			(D-\delta_D)([x_1, ..., x_n])=\e(d, X_i)[\a^k(x_1), ...,(D-\delta_D)(x_i), ...,\a^k(x_n) ],~i=1, ..., n.	
		\end{equation*}
		
		\item[(3)]  $\delta_{D-\delta_D}=-n(D-\delta_D).$	
	\end{itemize}	
\end{pro}	
\noindent {\bf Proof.} 
(1)  Suppose that $D\in\dd_{\a^k}(\LL).$ Since $\LL$ is perfect,  we can write  
$
x_j=\sum_{i\in{I_{j}}}[x_{i_{j_1}}, ..., x_{i_{j_n}}], j=1, ..., n-1$ and $y_s=\sum_{p\in{P_s}}[y_{p_{s_1}}, ..., y_{p_{s_n}}], s=1, ...,n.
$ 
From Eq. (\ref{122}), we have
$$
\delta_D(x_j)=\sum_{l=1}^{n}\sum_{i\in{I_{j}}}\e(d, X_{I_{j_l}})[\a^k(x_{i_{j_1}}), ..., D(x_{i_{j_l}}), ..., \a^k (x_{i_{j_n}})],	
$$
and similarly for $\delta_D(y_s).$ Next, we obtain
\begin{eqnarray*}
&\sum_{j=1}^{n-1}\e(d, X_j)\bigg[\a^k(x_1), ..., \delta_D(x_j), ..., \a^k(x_{n-1}), \big[\a^k(y_1), ..., \a^k(y_{n})\big]\bigg]\\
&+\sum_{s=1}^{n}\e(d, X+Y_s)\bigg[\a^k(x_1), ..., \a^k(x_{n-1}), \big[\a^k(y_1), ..., \delta_D(y_s), ..., \a^k(y_{n})\big]\bigg]\\
=&\sum_{j=1}^{n-1}\e(d, X_j)\bigg[\a^k(x_1), ..., \sum_{l=1}^{n}\sum_{i\in{I_{j}}}\e(d, X_{I_{j_l}})\\
&~[\a^k(x_{i_{j_1}}), ..., D(x_{i_{j_l}}), ..., \a^k (x_{i_{j_n}})], ..., \a^k(x_{n-1}), \big[\a^k(y_1), ..., \a^k(y_{n})\big]\bigg]\\
&+\sum_{s=1}^{n}\e(d, X+Y_s)\bigg[\a^k(x_1), ..., \a^k(x_{n-1}), \big[\a^k(y_1), ..., \sum_{l=1}^{n}\e(d, Y_{p_{s_l}})\\ 
&~[\a^k(y_{p_{s_1}}), ..., D(y_{p_{s_l}}), ..., \a^k(y_{n})], ..., \a^k(y_{n})\big]\bigg]\\
=&\sum_{i\in{I_{j}}}\bigg(\sum_{j=1}^{n-1}(-1)^{(n-1)(n-j)}\e(d, X_j)\e(Y, d+X)\e(\sum_{t=j+1}^{n-1}(x_t, \hat{X_t}+Y+d))\\
&.\bigg[\a^k(x_{j+1}), ..., \a^k(x_{n-1}), \big[\a^k(y_1), ..., \a^k(y_{n})], \a^k(x_1), ..., \a^k(x_{j-1}),\\
&~\sum_{l=1}^n\e(d, X_{i_{j_l}})[\a^k(x_{i_{j_1}}), ..., D(x_{i_{j_l}}), ..., \a^k (x_{i_{j_n}})]\bigg]\bigg)\\
&+\sum_{p\in{P_s}}\bigg(\sum_{s=1}^{n}(-1)^{(n-1)(n-j)}\e(d, X+Y_s)\e(\sum_{t=s+1}^{n}(y_t, \hat{Y_t}+d))\\ 
&.\bigg[\a^k(x_1), ..., \a^k(x_{n-1}), \big[\a^k(y_{s+1}), ..., \a^k(y_n), \a^k(y_1), ..., \a^k(y_{s_t}),\\
&~ \sum_{l=1}^{n}\e(d, Y_{p_{s_l}})[\a^k(y_{p_{s_1}}), ..., D(y_{p_{s_l}}), ..., \a^k(y_{p_{s_n}})]\big]\bigg]\bigg)\\
=&\sum_{i\in{I_{j}}}\bigg(\sum_{j=1}^{n}(-1)^{(n-1)(n-j)}\e(d, X_j)\e(Y, d+X)\\
&.\e(\sum_{t=j+1}^{n-1}(x_t, \hat{X_t}+Y+d))\e(d,\sum_{t=j+1}^{n-1}(\hat{X_t}+Y))\\
&.\bigg(D\big(\big[x_{j+1}, ..., x_{n-1}, [y_1, ..., y_{n}], x_1, ..., x_{j-1}, [x_{i_{j_1}}, ..., x_{i_{j_n}}]\big]\big)
\end{eqnarray*}
\begin{eqnarray*}
&-\sum_{m=j+1}^{n-1}\e(d, x_{j+1}+...+x_{m-1})\bigg[\a^k(x_{j+1}), ...,D(x_m),...,  \a^k(x_{n-1}),\\
&~~  [\a^k(y_1), ..., \a^k(y_n)], \a^k(x_1), ..., \a^k(x_{j-1}), [x_{i_{j_1}}, ..., x_{i_{j_n}}]\bigg]\\
&-\sum_{m=1}^{j-1}\e(d, x_{j+1}+...+x_{n-1}+Y+X_m)\bigg[\a^k(x_{j+1}), ..., \a^k(x_{n-1}),\\
&~ [\a^k(y_1), ..., \a^k(y_n)], \a^k(x_1), ...,D(x_m), ..., \a^k(x_{j-1}), [x_{i_{j_1}}, ..., x_{i_{j_n}}]\bigg]\\
&-\e(d, x_{j+1}+...+x_{n-1})\bigg[\a^k(x_{j+1}), ..., \a^k(x_{n-1}), D([y_1, ..., y_n]),\\
&~\a^k(x_1), ..., \a^k(x_{j-1}), [x_{i_{j_1}}, ..., x_{i_{j_n}}]\bigg]\bigg)\bigg)\\
&+\sum_{p\in{P_s}}\bigg(\sum_{s=1}^{n}(-1)^{(n-1)(n-j)}\e(d, X+Y_s)\\
&.\e(y_s, \sum_{t=s+1}^{n} \hat{Y_t}+d)\e(d, \sum_{t=s+1}^{n}\hat{Y_t})\\
&.\bigg(\big[\a^k(x_1), ..., \a^k(x_{n-1}), D\big([y_{s+1}, ..., y_n, y_1, ..., y_{s-1}, [y_{p_{s_1}}, ..., y_{p_{s_n}}]]\big)\big]\\
&-\sum_{m=s+1}^{n}\e(d, y_{s+1}+...+y_{m-1})\bigg[\a^k(x_1), ...,\a^k(x_{n-1}), \big[\a^k(y_{s+1}), ...,\\
&~ D(y_m), ..., \a^k(y_n), \a^k(y_1), ..., \a^k(y_{s-1}),[\a^k(y_{p_{s_1}}), ..., \a^k(y_{p_{s_n}})]\big]\bigg]\\
&-\sum_{m=1}^{s-1}\e(d, y_{s+1}+...+y_{n}+Y_m)\bigg[\a^k(x_1), ..., \a^k(x_{n-1}), [\a^k(y_{s+1}), ...,\\
&- \a^k(y_n), \a^k(y_1), ..., D(y_m), ...,  \a^k(y_{s-1}), [\a^k(y_{p_{s_1}}), ..., \a^k(y_{p_{s_n}})]\big]\bigg]\bigg)\bigg)\\
=&\sum_{j=1}^{n-1}(-1)^{(n-1)(n-j)}\e(d, X_j)\e(Y, d+X)\e(\sum_{t=j+1}^{n-1}(x_t, \hat{X_t}+Y+d))\\
&.\bigg((-1)^{(n-1)(n-j)}\e(\sum_{t=j+1}^{n-1}(x_t, \hat{X_t}+Y+d))\e(Y, X)\\
&.D\big(\big[x_1, ..., x_{n-1}, [y_1, ..., y_n]\big]\big)\\
&-\sum_{m=j+1}^{n-1}\e(d, x_{j+1}+...+x_{m-1})\e(d+x_m, \hat{X_m}+Y)\\
&.\e(\sum_{t=j+1}^{n-1}(x_t, \hat{X_t}+Y+d))\e(Y, d+X)\\
&.\big[\a^k(x_1), ..., \a^k(x_{n-1}), [\a^k(y_1), ..., \a^k(y_n)]\big]\\
&-\sum_{m=1}^{j-1}(-1)^{(n-1)(n-j)}\e(d, x_{j+1}+...+x_{n-1}+Y+X_m)\\
&.\e(\sum_{t=j+1}^{n-1}(x_t, \hat{X_t}+Y+d))\e(Y, X)
\end{eqnarray*}
\begin{eqnarray*}
&.\big[\a^k(x_1), ..., \a^k(x_{n-1}), [\a^k(y_1), ..., \a^k(y_n)]\big]\\
&-(-1)^{(n-1)(n-j)}\e(d, x_{j-1}+...+x_{n-1})\e(\sum_{t=j+1}^{n-1}(x_t, \hat{X_t}+Y+d))\\
&.\e(d+Y, X)\big[\a^k(y_1), ..., \a^k(y_n),  D\big([y_1, ..., y_n]\big)\big]\bigg)\\
&+\sum_{s=1}^{n}(-1)^{(n-1)(n-j)}\e(d, X+Y_s)\e(\sum_{t=j+1}^{n-1}(y_t, \hat{Y_t})\\
&.\bigg((-1)^{(n-1)(n-j)}\e(\sum_{t=j+1}^{n-1}(y_t, \hat{Y_t}))\\
&.\big[\a^k(y_1), ..., \a^k(y_n),  D\big([y_1, ..., y_s, ..., y_n]\big)\big]\\
&-\sum_{m=s+1}^{n}\e(d, y_{s+1}+...+y_{m-1})\e(\sum_{\substack{t=s+1\\ t\neq m}}^n(y_t, d+\hat{Y_t}))\\
&.\e(d+y_m, \hat{Y_m})\big[\a^k(y_1), ..., \a^k(y_n), [\a^k(y_1), ..., D(y_m), ..., \a^k(y_n)]\big]\\
&-\sum_{m=1}^{s-1}(-1)^{(n-1)(s-1)}\e(d, y_{s+1}+...+y_{n}+Y_m)\e(\sum_{t=s+1}^n(y_t, d+\hat{Y_t}))\\
&.\big[\a^k(x_1), ..., \a^k(x_{n-1}), [\a^k(y_1), ..., D(y_m), ..., \a^k(y_n)]\big]\bigg)\\
=&\sum_{j=1}^{n-1}D\bigg(\big[x_1, ...,x_j, ..., x_{n-1}, [y_1, ..., y_n]\big]\bigg)\\
&-\sum_{j=1}^{n-1}\sum_{m=j+1}^{n-1}\e(d, X_m)\\
&.\bigg[\a^k(x_1), ...,\a^k(x_j), ..., D(x_m), ..., \a^k(x_{n-1}), \big[\a^k(y_1), ..., \a^k(y_n)\big]\bigg]\\
&-\sum_{j=1}^{n-1}\sum_{m=1}^{j-1}\e(d, X_m)\\
&.\bigg[\a^k(x_1), ..., D(x_m), ..., \a^k(x_j), ..., \a^k(x_{n-1}), \big[\a^k(y_1), ..., \a^k(y_n)\big]\bigg]\\
&-\sum_{j=1}^{n-1}\e(d, X)\bigg[\a^k(x_1), ..., \a^k(x_j), ..., \a^k(x_{n-1}),  D\big([y_1, ...,  y_n]\big)\bigg]\\
&+\sum_{s=1}^{n}\e(d, X)\bigg[\a^k(x_1), ..., \a^k(x_{n-1}),  D\big([y_1, ..., y_s, ...,  y_n]\big)\bigg]\\
&-\sum_{s=1}^{n}\sum_{m=s+1}^{n}\e(d, X+Y_m)\\
&.\bigg[\a^k(x_1), ..., \a^k(x_{n-1}),  \big[\a^k(y_1), ..., \a^k(y_s), ...,  D(y_m), ..., \a^k(y_n)\big]\bigg]\\
&-\sum_{s=1}^{n}\sum_{m=1}^{s-1}\e(d, X+Y_m)\\
&.\bigg[\a^k(x_1), ..., \a^k(x_{n-1}),  \big[\a^k(y_1), ..., D(y_m), ..., \a^k(y_s), ..., \a^k(y_n)\big]\bigg]\\
=&\sum_{j=1}^{n-1}D\bigg(\big[x_1, ...,x_j, ..., x_{n-1}, [y_1, ..., y_n]]\big]\bigg)\\
&-\sum_{j=1}^{n-1}\sum_{\substack{m=1\\ m\neq j}}^{n-1}\e(d, X_m)\\
&.\bigg[\a^k(x_1), ...,D(x_m), ..., \a^k(x_{n-1}), \big[\a^k(y_1), ..., \a^k(y_n)\big]\bigg]
\end{eqnarray*}
\begin{eqnarray*}
&+\e(d, X)\bigg[\a^k(x_1), ..., \a^k(x_{n-1}), D\big([y_1, ..., y_n]\big)\bigg]\\
&-\sum_{s=1}^{n-1}\sum_{m=1}^{n}\e(d, X_m)\\
&.\bigg[\a^k(x_1), ..., \a^k(x_{n-1}), \big[\a^k(y_1), ..., D(y_m), ..., \a^k(y_n)\big]\bigg]\\
=&(n-1)D\bigg(\big[x_1, ...,x_j, ..., x_{n-1}, [y_1, ..., y_n]\big]\bigg)\\
&-(n-2)\sum_{m=1}^{n}\e(d, X_m)\\
&.\bigg[\a^k(x_1), ..., D(x_m), ..., \a^k(x_{n-1}), \big[\a^k(y_1), ..., \a^k(y_n)\big]\bigg]\\
&+\e(d, X)\bigg[\a^k(x_1), ..., \a^k(x_{n-1}), D\big([y_1, ..., y_n]\big)\bigg]\\
-&(n-1)\sum_{m=1}^{n-1}\e(d, X_m)\\
&.\bigg[\a^k(x_1), ..., D(x_m), ..., \a^k(x_{n-1}), \big[\a^k(y_1), ..., \a^k(y_n)\big]\bigg]\\
=&(n-1)\bigg(\sum_{m=1}^{n-1}\e(d, X_m)\\
&.\bigg[\a^k(x_1), ..., D(x_m), ..., \a^k(x_{n-1}), \big[\a^k(y_1), ..., \a^k(y_n)\big]\bigg]\\
&+\sum_{m=1}^{n}\e(d, X+Y_m)\\
&.\bigg[\a^k(x_1), ..., \a^k(x_{n-1}), \big[\a^k(y_1), ..., D(y_m), ..., \a^k(y_n)\big]\bigg]\bigg)\\
&-(n-2)\sum_{m=1}^{n-1}\e(d, X_m)\\
&.\bigg[\a^k(x_1), ..., D(x_m), ..., \a^k(x_{n-1}), \big[\a^k(y_1), ..., \a^k(y_n)\big]\bigg]\\
&+\e(d, X)\bigg[\a^k(x_1), ..., \a^k(x_{n-1}), D\big([y_1, ..., y_n]\big)\bigg]\\
&-(n-1)\sum_{m=1}^{n}\e(d, X+Y_m)\\
&.\bigg[\a^k(x_1), ..., \a^k(x_{n-1}), \big[\a^k(y_1), ..., D(y_m), ..., \a^k(y_n)\big]\bigg]\\
=&\sum_{m=1}^{n-1}\e(d, X_m)\big[\a^k(x_1), ..., D(x_m), ..., \a^k(x_{n-1}), [\a^k(y_1), ..., \a^k(y_n)]\big]\\
&+\e(d, X)\bigg[\a^k(x_1), ..., \a^k(x_{n-1}), D\big([y_1, ..., y_n]\big)\bigg]\\
=&\delta_D\bigg(\big[x_1, ..., x_{n-1}, [y_1, ..., y_n]\big]\bigg),
\end{eqnarray*}
where $ \hat{X_t} =\sum_{j=1}^{n-1} y_j.$ Thus, $\delta_D\in\dd_{\a^k}(\LL),$ and so is a double derivation of $\LL.$\\

(2)  For any $D\in\dd_{\a^k}(\LL)$ and $x_i, 
y_j, z_i\in\hh(\LL)$ with $1 \leq i\leq n - 1, 1\leq j \leq n,$ it follows
from Eq. (\ref{122}) that
\begin{eqnarray}
	&\nonumber\bigg[\a^k(z_1), ... , \a^k(z_{n-1}), D\big([x_1, ...,  x_{n-1}, [y_1, ...,  y_n]]\big)\bigg]\\
	\nonumber=&\sum_{i=1}^{n-1}\e(d, X_i)\bigg[\a^k(z_1), ... , \a^k(z_{n-1}),\\
\nonumber&~~\big[\a^k(x_1), ..., D(x_i), ...  \a^k(x_{n-1}), [\a^k(y_1), ..., \a^k(y_n)]\big]\bigg]\\
	\nonumber&+\sum_{j=1}^{n}\e(d, X+Y_j)\bigg[\a^k(z_1), ... , \a^k(z_{n-1}),\\
\nonumber&~\big[\a^k(x_1), ...  \a^k(x_{n-1}), [\a^k(y_1), ..., D(y_j), ..., \a^k(y_n)]\big]\bigg]\\
	\nonumber=&\bigg[\a^k(z_1), ... , \a^k(z_{n-1}), \delta_D\big([x_1, ...,  x_{n-1}, [y_1, ...,  y_n]]\big)\bigg]\\
	\nonumber&-\e(d, X)\bigg[\a^k(z_1), ... , \a^k(z_{n-1}), \big[\a^k(x_1), ...  \a^k(x_{n-1}), D\big([y_1, ..., y_n]\big)\big]\bigg]\\
	\nonumber&+\e(d, X)\bigg[\a^k(z_1), ... , \a^k(z_{n-1}), \big[\a^k(x_1), ...  \a^k(x_{n-1}), \delta_D\big([y_1, ..., y_n]\big)\big]\bigg].
\end{eqnarray}
Hence,
\begin{eqnarray}
	&\nonumber\bigg[\a^k(z_1), ... , \a^k(z_{n-1}), (D-\delta_D)\big([x_1, ...,  x_{n-1}, [y_1, ...,  y_n]]\big)\bigg]\\
	\nonumber=&-\bigg[\a^k(z_1), ... , \a^k(z_{n-1}), \big[\a^k(x_1), ...,  \a^k(x_{n-1}),\\
\nonumber&~~ \e(d, X)(D-\delta_D)([y_1, ...,  y_n])\big]\bigg].
\end{eqnarray}
Taking into account  $[\LL, \LL, ..., \LL]=\LL$ and $Z(\LL)=0$ we have that
\begin{eqnarray*}
	&(D-\delta_D)\big([x_1, ...,  x_{n-1}, [y_1, ...,  y_n]]\big)=\\
	&-\e(d, X)\big[\a^k(x_1), ...,  \a^k(x_{n-1}),  (D-\delta_D)\big([y_1, ...,  y_n]\big)\big],
\end{eqnarray*}
and so
 
\begin{equation*}
	(D-\delta_D)([x_1, ..., x_n])=\e(d, X_i)[\a^k(x_1), ...,(D-\delta_D)(x_i), ...,\a^k(x_n) ],	
\end{equation*}
for all $x_1, ..., x_n\in\hh(\LL)$ and $i=1, ..., n.$\\

(3) By part (1), if $D$ is
a double derivation of $\LL,$ then $\delta_D$ is a double derivation of $\LL.$ Hence, $D-\delta_D$ is
a double derivation of $\LL.$ Now, by part (2) we find that
\begin{eqnarray}
	\nonumber\delta_{D-\delta_D}([x_1, ...,  x_{n}])&=&\sum_{i=1}^n\e(d, X_i)[\a^k(x_1), ..., (D-\delta_D)(x_i), ..., \a^k(x_n)]\\
	\nonumber&=&-n(D-\delta_D)([x_1, ...,  x_{n}]),
\end{eqnarray}
which gives us the result.\qed

\begin{pro}\label{..} If $(\LL, [., ..., .], \a, \e)$ is a centerless multiplicative perfect $n$-Hom-Lie color algebra. Then the following assertions hold.
	\begin{itemize}
		\item[(1)]	 If $D$ is a double derivation of $\LL,$ then $\delta_D$ is a derivation of $\LL$ if and only if $D$
		is a derivation of $\LL.$ In particular, $\delta_D=D$ if $D$ is a derivation of $\LL.$
		\item[(2)] For all $x_1, ..., x_n\in\hh(\LL)$ and any $\a^k$-double derivation $D$ of $\LL,$ 
\begin{eqnarray*}
[D, \ad_s(x_1, ..., x_n)]&=&\ad_{k+s}(\delta_D(x_1), x_2, ..., x_{n-1})\\
&+&\sum_{j=2}^{n-1}\e(d, X_j)\ad_{k+s}(x_1,  ..., D(x_j), ..., x_{n-1}).
\end{eqnarray*}
	\end{itemize}	
\end{pro}	
\noindent {\bf Proof.} (1) Let $D$ be a derivation of $\LL.$ Then by Eq. (\ref{122}),
$$
\delta_D([x_1, ..., x_n])=D([x_1, ..., x_n])
$$
Taking into account $[\LL, \LL, ..., \LL]=\LL,$ we get $\delta_D=D.$ Next, suppose that $D\in\dd_{\a^k}(\LL)$ but not contained in $Der(\LL)$ and so $\delta_D\neq D.$ Then there exist $x_1, ..., x_n\in\hh(\LL)$ such that $(D -\delta_D)([x_1, ..., x_n])\neq 0.$ Suppose that $$x_j=\sum_{i\in{I_{j}}}[x_{i_{j_1}}, ..., x_{i_{j_n}}],~~~j=1, ..., n.$$  We have
\begin{eqnarray*}
&\sum_{j=1}^{n}\e(d, X_j)[\a^k(x_1), ..., \delta_D(x_j), ...  \a^k(x_{n})]\\
=&\sum_{j=1}^{n}\e(d, X_j)\bigg[\a^k(x_1), ..., \sum_{i\in{I_{j}}}\sum_{l=1}^{n}\e(d, X_{i_{j_l}})\big[\a^k(x_{i_{j_1}}),\\
~&~~ ..., D(x_{i_{j_l}}), ..., \a^k(x_{i_{j_n}})\big], ...  \a^k(x_{n})\bigg]\\
=&\sum_{j=1}^{n}\sum_{i\in{I_{j}}}(-1)^{(n-1)(n-j)}\e(d, X_j)\e(\sum_{m=j+1}^n(x_m, d+\hat{X_m}))\\
&.\bigg[\a^k(x_{j+1}), ..., \a^k(x_n), \a^k(x_1), ..., \sum_{l=1}^n\e(x, X_{i_{j_l}})\big[\a^k(x_{i_{j_1}}),\\
~&~~ ..., D(x_{i_{j_l}}), ..., \a^k(x_{i_{j_n}})\big]\bigg]\\
=&\sum_{j=1}^{n}\sum_{i\in{I_{j}}}(-1)^{(n-1)(n-j)}\e(d, X_j)\e(\sum_{m=j+1}^n(x_m, d+\hat{X_m}))\e(d, \hat{X_j})\\ 
&. \bigg(D\big([x_{j+1}, ..., x_n, x_1, ..., x_{j-1}, [x_{i_{j_1}}, ..., x_{i_{j_n}}]]\big)\\
&-\sum_{s=j+1}^{n}\e(d, x_{j+1}+...+x_{s-1})\\
&.[\a^k(x_{j+1}), ..., D(x_s), ..., \a^k(x_n), \a^k(x_1), ..., \a^k(x_{j-1}), [\a^k(x_{i_{j_1}}), ..., \a^k(x_{i_{j_n}})]]\\
&-\sum_{s=1}^{j-1}\e(d, x_{j+1}+...+x_{s-1})\\
&.[\a^k(x_{j+1}), ..., \a^k(x_n), \a^k(x_1), ..., D(x_s), ..., \a^k(x_{j-1}), [\a^k(x_{i_{j_1}}), ..., \a^k(x_{i_{j_n}})]]\bigg)\\
=&\sum_{j=1}^{n}(-1)^{(n-1)(n-j)}\e(d, X_j)\e(\sum_{m=j+1}^n(x_m, \hat{X_m}))\e(d, \hat{X_j})\\ 
&.\bigg((-1)^{(n-1)(n-j)}\e(\sum_{m=j+1}^n(x_m, \hat{X_m})) D\big([x_1, ..., x_n]\big)\\ 
&-\sum_{s=j+1}^{n}\e(d, x_{j+1}+...+x_{s-1})\e(\sum_{m=j+1}^n(x_m, d+\hat{X_m}))\e(d+x_s, \hat{ X_s})\\ 
&.[\a^k(x_1), ..., D(x_s), ..., \a^k(x_n)]\\ 
&-\sum_{s=1}^{j-1}(-1)^{(n-1)(n-j)}\e(d, x_{j+1}+...+x_{s-1})\e(\sum_{m=j+1}^n(x_m, d+\hat{X_m}))\\ 
&.[\a^k(x_1), ..., D(x_s), ..., \a^k(x_n)]\bigg)\\
=&\sum_{j=1}^{n}D\big([x_1, ..., x_n]\big)\\
&-\sum_{j=1}^{n}\sum_{s=j+1}^{n}\e(d, X_s)[\a^k(x_1), ..., D(x_s), ..., \a^k(x_n)]\\ 
&-\sum_{j=1}^{n}\sum_{s=1}^{j-1}\e(d, X_s)[\a^k(x_1), ..., D(x_s), ..., \a^k(x_n)]\\
=&nD\big([x_1, ..., x_n]\big)\\
&-\sum_{j=1}^{n}\sum_{{\substack{s=1\\ s\neq j}}}^{n}\e(d, X_s)[\a^k(x_1), ..., D(x_s), ..., \a^k(x_n)]\\
=&nD\big([x_1, ..., x_n]\big)-(n-1)\sum_{s=1}^{n}\e(d, X_s)[\a^k(x_1), ..., D(x_s), ..., \a^k(x_n)]\\
=&nD\big([x_1, ..., x_n]\big)-(n-1)\delta_D\big([x_1, ..., x_n]\big)\\
=&n(D-\delta_D)\big([x_1, ..., x_n]\big)+\delta_D\big([x_1, ..., x_n]\big)\\
\neq&\delta_D\big([x_1, ..., x_n]\big).
\end{eqnarray*}
It follows that $\delta_D\notin Der_{\a^k}(\LL)$ and so $\delta_D$ is not a derivation of $\LL.$\\

(2) For all $x_1, ..., x_n\in\hh(\LL)$ and any  $D\in\dd_{\a^k}(\LL),$  by Theorem \ref{ideal} and Eq. (\ref{122}), we have
\begin{eqnarray*}
&[D, ad_{s}(x_1, ...,  x_{n})]\\
=&\sum_{j=2}^{n-1}\e(d, X_j)\ad_{k+s}\big(x_1,  ..., D(x_j), ..., x_{n-1}\big)\\
&+\sum_{i\in I}\sum_{t=1}^{n}\e(d, X_{i_t})\ad_{k+s}\big([x_{i_1}, ..., D(x_{i_t}), ..., x_{i_n}], x_2, ..., x_{n-1}\big)\\
=&\sum_{j=2}^{n-1}\e(d, X_j)\ad_{k+s}\big(x_1,  ..., D(x_j), ..., x_{n-1}\big)\\
&+\sum_{i\in I}\sum_{t=1}^{n}\e(d, X_{i_t})\ad_{k+s}\big([\a^k(x_{i_1}), ..., D(x_{i_t}), ..., \a^k(x_{i_n})], x_2, ..., x_{n-1}\big)\\
=&\sum_{j=2}^{n-1}\e(d, X_j)ad_{k+s}\big(x_1,  ..., D(x_j), ..., x_{n-1}\big)+\ad_{k+s}\big(\delta_D(x_1), x_2, ..., x_{n-1}\big).
\end{eqnarray*}
Therefore, the result holds.\qed

\begin{pro} If $(\LL, [., ..., .], \a, \e)$ is a centerless perfect $n-$-Hom-Lie color algebra. Then the map $\delta ~:~\dd(\LL)\longrightarrow End(\LL)$ is a Hom-Lie color algebra homomorphism, that is, $\delta_{[D_1, D_2]}=[\delta_{D_1}, \delta_{D_2}].$
\end{pro}	
\noindent {\bf Proof.} Suppose that $D_1\in\dd_{\a^k}(\LL)$ and $D_2\in\dd_{\a^s}(\LL).$ Since $\LL$ is perfect, fory any $x\in\hh(\LL)$ we can write  $x=\sum_{i\in I}[x_{i_1}, ..., x_{i_n}]\in\LL.$ By Eq. (\ref{122}) we have
\begin{eqnarray*}
&\big[\delta_{D_1}, \delta_{D_2}\big](x)=\big(\delta_{D_1}\delta_{D_2}-\e(d_1, d_2)\delta_{D_2}\delta_{D_1}\big)(x)\\
=&\sum_{i\in I}\bigg(\delta_{D_1}\big(\delta_{D_2}([x_{i_1}, ..., x_{i_n}])\big)-\e(d_1, d_2)\delta_{D_2}\big(\delta_{D_1}([x_{i_1}, ..., x_{i_n}])\big)\bigg)\\
=&\sum_{i\in I}\bigg(\delta_{D_1}\big(\sum_{j=1}^{n}\e(d_2, X_{i_j})[\a^s(x_{i_1}), ..., D_2(x_{i_j}), ..., \a^s(x_{i_n})]\big)\\ 
&-\e(d_1, d_2)\delta_{D_2}\big(\sum_{j=1}^{n}\e(d_1, X_{i_j})[\a^k(x_{i_1}), ..., D_1(x_{i_j}), ..., \a^k(x_{i_n})]\big)\bigg)\\
=&\sum_{i\in I}\bigg(\sum_{j=1}^{n}\big(\sum_{l=1}^{j-1}\e(d_2, X_{i_j})\e(d_1, X_{i_l})\\
&.\big[\a^{k+s}(x_{i_1}), ..., D_1(\a^k(x_{i_l})), ...,  D_2(x_{i_j}), ..., \a^{k+s}(x_{i_n})\big]\\
&+\sum_{l=j+1}^{n}\e(d_2, X_{i_j})\e(d_1, X_{i_l})\\
&.[\a^{k+s}(x_{i_1}), ..., D_2(\a^s(x_{i_l})), ...,  D_1(x_{i_j}), ..., \a^{k+s}(x_{i_n})]\big)
\end{eqnarray*}
\begin{eqnarray*}
&-\e(d_1, d_2)\sum_{j=1}^{n}\big(\sum_{l=1}^{j-1}\e(d_1, X_{i_j})\e(d_2, X_{i_l})\\
&.[\a^{k+s}(x_{i_1}), ...,  D_2(\a^s(x_{i_l})), ..., D_1(x_{i_j}), ..., \a^{k+s}(x_{i_n})]\\
&+\sum_{l=j+1}^{n}\e(d_2, X_{i_j})\e(d_1, X_{i_l})\\
&.[\a^{k+s}(x_{i_1}), ..., D_1(\a^k(x_{i_l})), ...,  D_2(x_{i_j}), ..., \a^{k+s}(x_{i_n})]\bigg)\\
=&\sum_{i\in I}\bigg(\sum_{j=1}^{n}\e(d_1+d_2, X_{i_j})[\a^{k+s}(x_{i_1}), ..., D_1D_2(x_{i_j}), ..., \a^{k+s}(x_{i_n})]\\
&-\sum_{j=1}^{n}\e(d_1+d_2, X_{i_j})[\a^{k+s}(x_{i_1}), ..., \e(d_1, d_2)D_2D_1(x_{i_j}), ..., \a^{k+s}(x_{i_n})]\bigg)\\
=&\sum_{i\in I}\sum_{j=1}^{n}\e(d_1+d_2, X_{i_j})\\
&.\bigg[\a^{k+s}(x_{i_1}), ..., \big(D_1D_2-\e(d_1, d_2)D_2D_1\big)(x_{i_j}), ..., \a^{k+s}(x_{i_n})\bigg]\\
=&\sum_{i\in I}\sum_{j=1}^{n}\e(d_1+d_2, X_{i_j})\\
&.\bigg[\a^{k+s}(x_{i_1}), ..., \big[D_1, D_2\big](x_{i_j}), ..., \a^{k+s}(x_{i_n})\bigg]\\
=&\delta_{[D_1, D_2]}(x).
\end{eqnarray*}
\qed

\begin{thm}\label{cent} If $(\LL, [., ..., .], \a, \e)$ is a  perfect $n$-Hom-Lie color algebra. Then the centralizer of $Inn(\LL)$ in $\dd(\LL)$ is trivial. In particular, the center of the Hom-Lie color algebra $\dd(\LL)$ is zero.
\end{thm}	
\noindent {\bf Proof.} Let $D\in\dd_{\a^s}(\LL)$ of degree $d$ and $D\in C_{\dd(\LL)}(Inn(\LL)).$ Then 
$$
[D, ad_{k}(x_1, ...,  x_{n-1})]=0,~~\forall x_1, x_2, ..., x_{n-1}\in\hh(\LL).
$$
Now, for all $z\in\hh(\LL),$ we have 
\begin{eqnarray}
	\nonumber 0&=&[D, ad_{k}(x_1, ...,  x_{n-1})](z)\\
	\nonumber&=&D\circ ad_k(x_1, ...,  x_{n-1})(z)-\e(d, X)ad_k(x_1, ...,  x_{n-1})\circ D(z)\\
	\nonumber&=&D\big([x_1, ...,  x_{n-1}, \a^k(z)]\big)-\e(d, X)[[\a^k(x_1), ...,  \a^k(x_{n-1}), \a^k(D(z))].	
\end{eqnarray}
From here,  we obtain
\begin{eqnarray}\label{126}
	D\big([x_1, ...,  x_{n-1}, \a^k(z)]\big)&=&\e(d, X)[\a^k(x_1), ...,  \a^k(x_{n-1}), \a^k(D(z))]\\
	\nonumber	&=&\e(d, X_i)[\a^k(x_1), ..., D(x_i), ...,  \a^k(x_{n-1}), \a^k(z)],	
\end{eqnarray}
for $1\leq i\leq n-1.$

Next, let $x_1, ..., x_{n-1}, y_1, ..., y_n\in\hh(\LL)$ and  taking into account Eq. (\ref{126}), we get
\begin{eqnarray}
	\nonumber&	D\big([x_1, ..., x_{n-1}, [y_1, ..., y_{n}]]\big)\\
	\nonumber=&\sum_{i=1}^{n-1}\e(d, X_i)\big[\a^k(x_1), ..., D(x_i), ..., \a^k(x_{n-1}), [\a^k(y_1), ..., \a^k(y_n)]\big]\\
	\nonumber&+\sum_{j=1}^{n}\e(d, X+Y_{j})\big[\a^k(x_1), ..., \a^k(x_{n-1}), [\a^k(y_1), ..., D(y_j), ..., \a^k(y_n)]\big]\\
	\nonumber=&(n-1)	D\big([x_1, ..., x_{n-1}, [y_1, ..., y_{n}]]\big)+n D([x_1, ..., x_{n-1}, [y_1, ..., y_{n}]]\big)\\
	\nonumber=&(2n-1) D\big([x_1, ..., x_{n-1}, [y_1, ..., y_{n}]]\big),
\end{eqnarray}
so we obtain $(2n-2) D\big([x_1, ..., x_{n-1}, [y_1, ..., y_{n}]]\big)=0.$  That is,
$$
D\big([x_1, ..., x_{n-1}, [y_1, ..., y_{n}]]\big)=0.
$$
Taking into account $[\LL, \LL, ..., \LL]=\LL,$ we conclude that $D=0.$\qed

\section{Triple derivations of double derivation algebra $\dd(\LL)$}\setcounter{equation}{0}\

In the final section, we are particularly interested in the  triple derivations of the Hom-Lie color algebra
$(\dd(\LL), [., .], \ta, \e)$ consisting of double derivations algebra as well as in the triple derivations of
the Hom-Lie color algebra $(Inn(\LL), [., .], \ta, \e)$  consisting of inner derivations algebra of a perfect $n$-Hom-Lie color algebra $\LL.$

Recall that, a linear map $D$ of a multiplicative Hom-Lie color algebra $(\L, [., .], \a, \e)$ is called an $\a^k$-triple derivation if it satisfies\\

(1)  $D\circ\a=\a\circ D.$ 

(2) For all $x, y, z\in\hh(\L),$
\begin{eqnarray}\label{triple}	
	D\big([x, [y, z]]\big)&=& \big[D(x), [\a^k(y), \a^k(z)]\big]+\e(d, x) \big[\a^k(x), [D(y),  \a^k(z)]\big]\\
\nonumber&+&\e(d, x+y) \big[\a^k(x), [\a^k(y), D(z)]\big].
\end{eqnarray}	
Let us denote the set of all $\a^k$-triple derivations of  $\L$ by $TDer_{\a^k}(\L)$ and denote  $TDer(\L) = \bigoplus_{k\geq 0}TDer_{\a^k} (\L),$ which is a Hom-Lie color algebra and it is called the  triple derivation algebra of $\L$

In the following we suppose that $\LL$ is a centerless  perfect $n$-Hom-Lie color algebra and  $\L=(Inn(\LL), [., .], \ta=D\circ \a, \e)$ is a perfect Hom-Lie color algebra. Then we have the following results.

\begin{lem}\label{3.3.} If $(\LL, [., ..., .], \a, \e)$ is a  perfect $n$-Hom-Lie color algebra. Then every  triple derivation $D$ of $(\dd(\LL), [., .], \ta, \e)$  keeps $\L=Inn(\LL)$ invariant. Furthermore, if $D(\L) = 0,$ then $D = 0.$
\end{lem}	
\noindent {\bf Proof.} First, we show that $D(\L)\subseteq \L,$ for $D\in TDer_{\ta^k}(\dd(\LL)).$  Since $\L=\{\ad_s(x_1, ..., x_{n-1})~:~\a(x_i)=x_i\in\LL, i=1, ..., n-1\}$ is  a perfect Hom-Lie color algebra, for all $x_1, ..., x_{n-1}\in\hh(\LL)$  there exist $x_{1_i}, ..., x_{{n-1}_i}, y_{1_{i_{j}}}, ..., y_{{n-1}_{i_{j}}}\in\LL$ and $z_{1_{i_{j}}}, ..., z_{{n-1}_{i_{j}}}\in\hh(\LL)$ such that
\begin{eqnarray}	
	\nonumber& \ad_s(x_1, ..., x_{n-1})\\
	\nonumber=& \sum_{i\in I}[\ad_s(x_{1_i}, ..., x_{{n-1}_i}), \ad_s(y_{1_i}, ..., y_{{n-1}_i})]\\
	\nonumber=& \sum_{i\in I}\bigg(\big[\ad_s(x_{1_i}, ..., x_{{n-1}_i}), \sum_{j\in J}[\ad_s( y_{1_{i_{j}}}, ..., y_{{n-1}_{i_{j}}}), \ad_s(z_{1_{i_{j}}}, ..., z_{{n-1}_{i_{j}}})]\big]\bigg),
\end{eqnarray}
for some finite index sets $I$ and $J.$  Now, by definition of triple derivation we have
\begin{eqnarray}
	\nonumber&D\big(\ad_s(x_1, ..., x_{n-1})\big)\\	
	\nonumber& \sum_{i\in I}\sum_{j\in J}D\bigg(\big[\ad_s(x_{1_i}, ..., x_{{n-1}_i}), [ad_s( y_{1_{i_{j}}}, ..., y_{{n-1}_{i_{j}}}), \ad_s(z_{1_{i_{j}}}, ..., z_{{n-1}_{i_{j}}})]\big]\bigg)\\
	\nonumber=&\sum_{i\in I}\sum_{j\in J}\bigg(\big[D(\ad_s(x_{1_i}, ..., x_{{n-1}_i})),\\ \nonumber& [\ad_s(\ta^k( y_{1_{i_{j}}}), ..., \ta^k(y_{{n-1}_{i_{j}}})), \ad_s(\ta^k(z_{1_{i_{j}}}), ..., \ta^k(z_{{n-1}_{i_{j}}}))]\big]\\
	\nonumber&+\e(d, X_i) \big[\ad_s(\ta^k(x_{1_i}), ..., \ta^k(x_{{n-1}_i})),\\ \nonumber& [D(\ad_s(\ta^k(y_{1_{i_{j}}}), ..., \ta^k(y_{{n-1}_{i_{j}}}))), \ad_s(\ta^k(z_{1_{i_{j}}}), ..., \ta^k(z_{{n-1}_{i_{j}}}))]\big]\\
	\nonumber&+\e(d, X_i+Y_{i_{j}})\big[\ad_s(\ta^k(x_{1_i}), ..., \ta^k(x_{{n-1}_i})),\\ \nonumber& [\ad_s(\ta^k(y_{1_{i_{j}}}), ..., \ta^k(y_{{n-1}_{i_{j}}})),  D(\ad_s(z_{1_{i_{j}}}, ..., z_{{n-1}_{i_{j}}}))]\big]\bigg).
\end{eqnarray}
Thanks to Theorem \ref{ideal}, we obtain $D(ad_s(x_1, ..., x_{n-1}))$ is an inner $\ta^{k+s}$-derivation and so $D(\L)\subseteq \L.$

Next,  let $D$ be an $\ta^k$-triple derivation of $\dd(\LL)$ which satisfies $D(\L) = 0.$ Then for each double derivation $d$ of $\LL,$ again by perfectness of $\L$ we have
\begin{eqnarray}	
	\nonumber& [D(d), \ad_s(x_1, ..., x_{n-1})]\\
	\nonumber=& \sum_{i\in I}\big[D(d), [\ad_s(x_{1_i}, ..., x_{{n-1}_i}), \ad_s(y_{1_i}, ..., y_{{n-1}_i})]\big]\\
	\nonumber=& \sum_{i\in I}\bigg(D\big([d, [\ad_s(x_{1_i}, ..., x_{{n-1}_i}), \ad_s(y_{1_i}, ..., y_{{n-1}_i})]]\big)\\
	\nonumber&-\e(d, e)\big[\ta^k(d), [D\big(\ad_s(x_{1_i}, ..., x_{{n-1}_i})\big), \ta^k(\ad_s(y_{1_i}, ..., y_{{n-1}_i})]\big]\\
	\nonumber&-\e(d, e+X_i)\big[\ta^k(d), [\ta^k(\ad_s(x_{1_i}, ..., x_{{n-1}_i})), D\big(\ad_s(y_{1_i}, ..., y_{{n-1}_i})\big)]\big]\bigg).
\end{eqnarray}
Since $[d, \L]\subseteq \L$ (tanks to Theorem \ref{ideal}), we get $D\big([d, \ad_s(x_{1_i}, ..., x_{{n-1}_i})]\big) = 0.$ Thus, we can obtain
$[D(d), \ad_s(x_{1_i}, ..., x_{{n-1}_i})] = 0.$ Therefore, $D(d)$ is in the centralizer of $\L$ in $\dd(\LL),$
and by Proposition \ref{cent}, $D(d) = 0,~\forall d\in\dd(\LL).$ Hence, $D = 0.$ The proof is completed.\qed

\begin{lem}\label{3.3} Let $(\LL, [., ..., .], \a, \e)$ be a centerles perfect $n$-Hom-Lie color algebra. If $D \in TDer(Der(\LL)),$ then there exists $\Delta \in Der(\LL)$ such that $D(\ad_s(x))=\ad_s(\Delta(x))$ for all $x\in\hh(\LL).$
\end{lem}	
\noindent {\bf Proof.} For any $D\in TDer_{\ta^k}(Der(\LL)),$ by Lemma \ref{3.3.}, we can obtain  $D(ad_s(x))\in Inn_{\ta^{k+s}}(\LL)$ for all $x\in\hh(\LL).$ Now, let $y\in\hh(\LL)$ such that $D(\ad_s(x))=ad_{k+s}(y).$ Note that $y$ is  unique, thanks to $Z(\LL)=0.$ Define the map $\Delta~:~\LL\longrightarrow\LL$ by $\Delta(x)=y$ which is an endomorphism. For $x_1, x_2, x_3\in\hh(\LL),$ we have
\begin{eqnarray}	
	\nonumber& \ad_{k+s}\bigg(\Delta\big([[x_1, x_2], x_3]\big)\bigg)=D\bigg(\ad_s\big([[x_1, x_2], x_3]\big)\bigg)\\
	\nonumber=&D\bigg(\big[[\ad_s(x_1), \ad_(x_2)], \ad_s(x_3)\big]\bigg)\\
	\nonumber=& \big[[D(\ad_s(x_1)), \ta^k(\ad_(x_2))], \ta^k(\ad_s(x_3))\big]\\ 
\nonumber&+\e(d, x_1)\big[[\ta^k(\ad_s(x_1)) , D(\ad_(x_2))], \ta^k(\ad_s(x_3))\big]\\
\nonumber&+\e(d, x_1+x_2)\big[[\ta^k(\ad_s(x_1)), \ta^k(\ad_(x_2))], D(\ad_s(x_3))\big]\\
	\nonumber=& \big[[\ad_{k+s}(\Delta(x_1)), \ta^k(\ad_(x_2))], \ta^k(\ad_s(x_3))]\\ \nonumber&+\e(d, x_1)\big[[\ta^k(\ad_s(x_1)), \ad_{k+s}(\Delta(x_2))], \ta^k(\ad_s(x_3))\big]\\ \nonumber&+\e(d, x_1+x_2)\big[[\ta^k(\ad_s(x_1)), \ta^k(\ad_(x_2))], \ad_{k+s}(\Delta(x_3))\big]\\
	\nonumber=&\ad_{k+s}\bigg(\big[\Delta(x_1), [\a^k(x_2), \a^k(x_3)]\big]+\e(d, x_1) \big[\a^k(x_1), [\Delta(x_2),  \a^k(x_3)]\big]\\ \nonumber&+\e(d, x_1+x_2) \big[\a^k(x_1), [\a^k(x_2), \Delta(x_3)]\big]\bigg).
\end{eqnarray}
Since $Z(\LL) = 0,$ we obtain
\begin{eqnarray}\label{triple}	
	&\nonumber\Delta\bigg(\big[[x_1, x_2], x_3\big]\bigg)\\
\nonumber=& \big[[\Delta(x_1), \a^k(x_2)], \a^k(x_3)\big]+\e(d, x_1) \big[[\a^k(x_1), \Delta(x_2)],  \a^k(x_3)\big]\\ 
\nonumber&+\e(d, x_1+x_2) \big[[\a^k(x_1), \a^k(x_2)], \Delta(x_3)\big].
\end{eqnarray}
That is, $ \Delta\in TDer_{\ta^k}(\LL).$ By Proposition \ref{..}, $\Delta\in Der(\LL).$\qed

\begin{thm}\label{3.4} Let $(\LL, [., ..., .], \a, \e)$ be a centerles perfect $n$-Hom-Lie color algebra. Then the triple derivation	algebra of $Inn(\LL)$ coincides with the derivation algebra of $Inn(\LL).$
\end{thm}	
\noindent {\bf Proof.} We have  $\L=(Inn(\LL), [., .], \ta, \e)$ is a perfect Hom-Lie color algebra with zero center.  For any $D\in TDer_{\ta^k}(\L),$ by Lemma \ref{3.3} there exists $\Delta\in Der(\LL)$ such that  $D(ad_s(x))=ad_s(\Delta(x))$ for all $x\in\LL.$ As well known that  $ad_s(\Delta(x))=[\Delta, ad_s(x)].$ We also have
$$
D(ad_s(x))=ad_s(\Delta(x))=[\Delta, ad_s(x)]=ad_s(\Delta)(ad_s(x))
$$
Hence, $(D-ad_s(\Delta(x)))(ad_s(x))=0.$ By Lemma \ref{3.3.}, $D=ad_s(\Delta(x)).$ Therefore, $TDer(\L)$ equal to derivation algebra of $Inn(\LL).$\qed

\begin{thm} Let $(\LL, [., ..., .], \a, \e)$ be a  perfect $n$-Hom-Lie color algebra. Then the triple derivation algebra of Hom-Lie color algebra $(Der(\LL), [., .], \ta, \e)$ is equal to the derivation algebra  $Der(\LL).$
\end{thm}
\noindent {\bf Proof.} Suppose that  $D$ is a $\ta^k$-triple derivation of $Der(\LL),$ then $D$ is an $\ta^k$-triple derivation of $\L=Inn(\LL)$
by Lemma \ref{3.3}, we get  $D$ is an $\ta^k$-derivation of $\L=Inn(\LL),$ thanks to  Theorem \ref{3.4}.  Since $\L$ is a perfect Hom-Lie color algebra,  for any
$x_1, ..., x_{n-1}\in\hh(\LL),$ there exist $x_{1_i}, ..., x_{{n-1}_i}, y_{1_i}, ..., y_{{n-1}_i}\in\hh(\LL)$  such that
$$
\ad_s(x_1, ..., x_{n-1})=\sum_{i\in I}[\ad_s(x_{1_i}, ..., x_{{n-1}_i}), \ad_s(y_{1_i}, ..., y_{{n-1}_i})],
$$
for some finite index set $I.$ Then for all $d_1, d_2\in Der(\LL)$ of degrees $\lam_1, \lam_2$ respectively, we have
\begin{eqnarray}	
	\nonumber& D\bigg(\big[[d_1, d_2], \ad_s(x_1, ..., x_{n-1})\big]\bigg)\\
	\nonumber=&D\bigg(\big[[d_1, d_2],\sum_{i\in I} [\ad_s(x_{1_i}, ..., x_{{n-1}_i}), \ad_s(y_{1_i}, ..., y_{{n-1}_i})]\big]\bigg)\\
	\nonumber=&\bigg[D([d_1, d_2]), \sum_{i\in I} \big[\ad_s(\ta^k(x_{1_i}), ..., \ta^k(x_{{n-1}_i})), \ad_s(\ta^k(y_{1_i}), ..., \ta^k(y_{{n-1}_i}))\big]\bigg]\\
	\nonumber&-\e(d, \lam_1+\lam_2)\\
\nonumber&.\bigg[\ta^k([d_1, d_2]), \sum_{i\in I}\big[D(\ad_s(x_{1_i}, ..., x_{{n-1}_i})), \ad_s(\ta^k(y_{1_i}), ..., \ta^k(y_{{n-1}_i}))\big]\bigg]\\
	\nonumber&-\e(d, \lam_1+\lam_2+X_i)\\
\nonumber&.\bigg[\ta^k([d_1, d_2]), \sum_{i\in I}\big[\ad_s(\ta^k(x_{1_i}), ..., \ta^k(x_{{n-1}_i}))), D(\ad_s(y_{1_i}, ...,  y_{{n-1}_i}))\big]\bigg]\\
\nonumber=&\bigg[D\big([d_1, d_2]\big), \ta^k(\ad_s(x_1, ..., x_{n-1}))\bigg]\\
\nonumber&+\e(d, \lam_1+\lam_2)\bigg[\ta^k[d_1, d_2], D\big(\ad_s(x_1, ..., x_{n-1})\big)\bigg]\\
	\nonumber=&\bigg[D\big([d_1, d_2]\big), \ta^k(\ad_s(x_1, ..., x_{n-1}))\bigg]\\
\nonumber&-\e( \lam_1+\lam_2, X)D\bigg(\big[\ad_s(x_1, ..., x_{n-1}), [d_1, d_2]\big]\bigg)\\
	\nonumber&+\e(d, X)\e( \lam_1+\lam_2, X)\bigg[\ta^k(\ad_s(x_1, ..., x_{n-1})), \big[D(d_1), \ta^k(d_2)\big]\bigg]\\
	\nonumber&+\e(d, \lam_1+X)\e( \lam_1+\lam_2, X)\bigg[\ta^k(\ad_s(x_1, ..., x_{n-1})), \big[\ta^k(d_1), D(d_2)\big]\bigg]
\end{eqnarray}
\begin{eqnarray}
	\nonumber=&\bigg[D\big([d_1, d_2]\big), \ta^k(\ad_s(x_1, ..., x_{n-1}))\bigg]+D\bigg(\big[[d_1, d_2], \ad_s(x_1, ..., x_{n-1})\big]\bigg)\\
	\nonumber=&\bigg[\big[D(d_1), \ta^k(d_2)\big], \ta^k(\ad_s(x_1, ..., x_{n-1}))\bigg]\\
\nonumber&-\e(d, \lam_1)\bigg[\big[\ta^k(d_1), D(d_2)\big], \ta^k(\ad_s(x_1, ..., x_{n-1}))\bigg]. 
\end{eqnarray}
Thus,
$$
D\big([d_1, d_2]\big)-\big[D(d_1), \ta^k(d_2)\big]-\e(d, \lam_1)\big[[\ta^k(d_1), D(d_2)], \ta^k( ad_s(x_1, ..., x_{n-1}))\big]=0.
$$
Now, Proposition \ref{cent} gives us 
$$
D\big([d_1, d_2]\big)=[D(d_1), \ta^k(d_2)]-\e(d, \lam_1)[\ta^k(d_1), D(d_2)].
$$
That is $D$ is an $\ta^k$-derivation of Hom-Lie color algebra $Der(\LL).$\qed






\end{document}